 %%%%%%%%%%%%%%%%%%%%%%%%%%%%%%%%%%%%%%%%%%%%%%%%%%%%%%%%
 %
 % Use standard AMSTeX command \proclaim for Theorems, Lemmas and other
 % statements. Statements should be numbered consecutively in each section.
 %
 % Insert all your macros as indicated below.
 %
 % Note that, when referring to a particular section, chapter, theorem, etc., 
 % in the text, the word in question should be in capitals. For example:
 % In this section we prove Theorem 2.1 formulated in Section 2.
 %
 % All numbers and brackets {}, (), [] should be roman(!!!!!), even within
 % italic text.
 %
 % Fractions in displayed formulas: displaystyle size;
 % in text: textstyle size
 % 
 % The commands \rightheadtext and \leftheadtext may not work with older 
 % versions of AMS-TeX. 
 % 
 % References should be given a number rather than initials. 
 %
 %%%%%%%%%%%%%%%%%%%%%%%%%%%%%%%%%%%%%%%%%%%%%%%%%%%%%%%%
\input amstex
\documentstyle{amsppt}

\hsize=4.75in
\vsize=8in
\NoBlackBoxes
 %%%%%%%%%%%%%%%%%%%%%%%%%%%%%%
 %% insert your macros here %%
\def\Pg{{\Cal P}^{\uparrow}_+}
\def\Lx{\Lambda,x}
\def\Aa{{\Cal A}}
\def\Oo{{\Cal O}}
\def\Rl{{\Bbb R}}
\def\ax{\alpha_x}
\def\aLx{\alpha_{\Lambda,x}}
\def\aP{\alpha_{{\Cal P}^{\uparrow}_+}
\def\Aul{{\underline{A}_{\lambda}}}
\def\Alu{\underline{\Cal A}}
\def\Au{\underline{A}}
\def\aLxu{\underline{\alpha}_{\Lambda,x}}
\def\oul{\underline{\omega}_{\lambda}}
\def\opul{\underline{\omega}_{\lambda}'}
\def\lauP{{}^{(\lambda)}\!\underline{\alpha}_{{\Cal P}^{\uparrow}_+}}}
\def\Hh{{\Cal H}}
\def\Ff{{\Cal F}}

\def\Aoi{{\Cal A}_{0,\iota}}

 %%%%%%%%%%%%%%%%%%%%%%%%%%%%%%

\rightheadtext {Scaling Algebras}
\leftheadtext {Rainer Verch}
\topmatter
\title
Scaling Algebras, the Renormalization Group and the Principle of
Local Stability in Algebraic Quantum Field Theory
\endtitle
\author
Rainer Verch
\endauthor
\affil
Dipartimento di Matematica,\\  Universit\`a di Roma II ``Tor Vergata'',\\
I - 00133 Roma, Italy
\endaffil
\endtopmatter
\document
%%%
\heading
1. An Invitation to the Scaling Algebra
\endheading
%%%%%%%%%%%%%%%%%%%%%%%%%%%%%%%%%%%%%%%%%%%%%%%%%%%%%%%%%%%%%%%%%%%%%%%%%%
In the present contribution we report on a new approach to the 
structural analysis of the short distance behaviour of quantum field
theories in the operator-algebraic formulation which has recently
been proposed in [1] (see also [2,3,4]), and the extension of this
framework to quantum field theory in curved spacetimes as set out
in [5].
\par
I shall begin by giving a few indications as to why one is interested
in having a model-independent framework for the analysis of the
short-distance properties of general quantum field theories.
By general agreement, quantum field theory is so far the best
description of elementary particle physics available. In elementary
particle physics there is experimental evidence that the hadrons
(e.g., protons and neutrons) are built up from particle-like
constituents, the quarks and gluons which, however, are never
observed as free particles but are ``confined'' in their bound states.
 Only when hadrons participate in
collision processes at very high energies/short distances in
particle accelerators, traces of those particle-like sub-structures
become indirectly visible. In the --- up to now mainly perturbative ---
quantum field theoretical treatment  of
strong interaction processes (quantum chromodynamics) the phenomenon of
confinement, and the related notion of asymptotic freedom, can be
explained by the method of renormalization group analysis of the
short distance behaviour in the corresponding quantum field models.
(Standard references include e.g.\ [6].)

In its conventional form, the renormalization group analysis
{depends} on the use of Wightman-type quantum fields. This is quite
unsatisfactory from the point of view of structural analysis since the
description of a theory in terms of quantum fields is not intrinsic.
The field operators (operator-valued distributions) merely serve as
a ``coordinatization'' of the field- or observable-algebras, and there
is some arbitrariness in their choice [7]. Therefore, to have
e.g.\ an intrinsic criterion as to whether a given quantum field
theory obeys a dynamics where confinement occurs, one needs a formulation
of renormalization group analysis in local algebraic quantum field
theory which uses only the intrinsic notions of this setting, like
observables, states, localization properties, automorphism groups,
 superselection sectors, etc.

The same applies to two other important cases in which the short-distance
analysis by the renormalization-group method in terms of quantum
fields has led to insights of considerable conceptual and mathematical
interest. The first is the connection, established by Fredenhagen [8],
between a stability property of the short distance behaviour of a
Wightman field theory and the property of its local von Neumann
algebras to be of type III${}_1$. The second arises in the
context of quantum field theory in curved spacetime where, due to
the absence of spacetime symmetries in general, one faces difficulties
to fix the folium of physical states in a manner which sufficiently
captures the idea of dynamical stability underlying the spectrum
condition. One possibility is to demand that, the closer one gets
to an arbitrary point in spacetime, the more similar a theory should
become to a quantum field theory on Minkowski-spacetime, satisfying the
usual spectrum condition. (For another possibility, where the dynamical
stability requirement is formulated in terms of a ``microlocal spectrum
condition'', the reader is referred to the contribution by K. Fredenhagen
in this volume.) So invariance and spectrum condition is, for
physical states, asymptotically realized in the short-distance
``scaling limit''. This condition was introduced in [9] and called
``principle of local stability'' (see also [10,11]). One
of its interesting consequences is that it allows to fix the Hawking
temperature of a quantum field state in a black hole spacetime [9].

Let us look in some more detail at the renormalization group
transformations in terms of quantum fields, and at the resulting
short-distance scaling limits. At the very beginning we can be 
quite general and consider a quantum field over a curved spacetime.
To this end we recall that a spacetime $(M,g)$ consists of a
four-dimensional $C^{\infty}$-manifold $M$ together with a Lorentzian metric
$g$ on $M$. That means, $g$ is a smooth section in $T^*M \otimes T^*M$,
with the property that around each point $p \in M$ there are
coordinates $(x^{\mu})$ so that $(g_{\mu\nu} |_p) = \text{diag}(1,
-1,-1,-1)$. Such coordinates are called {\it normal coordinates}
at $p$ if they also map $p$ to 0, $x^{\mu}(p) = 0$.
 We will assume that $(M,g)$ is space- and time-orientable and that
such orientations have been chosen. Now consider a quantum field
$C_0^{\infty}(M) \owns h \mapsto \phi(h)$ over $(M,g)$, i.e.\ an
operator valued distribution defined on the test-functions over $M$
and taking values in the essentially selfadjoint operators having a
common dense and invariant domain in some Hilbert-space $\Cal H$.
Fix some point $p \in M$ and a normal coordinate system $(x^{\mu})$
around $p$. Then the ``standard type'' of renormalization group
transformations (with respect to these normal coordinates) is 
a family $R_{\lambda}$, $\lambda > 0$, on the field operators,
acting as follows:
$$ R_{\lambda} : \phi(f) \mapsto \phi_{\lambda}(f) :=
 N_{\lambda} \cdot \phi(f_{\lambda})\,,\quad \lambda > 0\,, \eqno(1.1) $$
where $f_{\lambda}(x) := f(\lambda^{-1}x)$ in the chosen normal
coordinate system, for each test-function $f$ supported sufficiently
close to $p$, and $N_{\lambda}$ is a positive real number depending
on $\lambda$. Thus $R_{\lambda}$ acts through scaling the spacetime
coordinates by the scaling factor $\lambda$, and by a ``multiplicative
field strength renormalization'', provided by the numerical factor
$N_{\lambda}$. Then look at the state $\langle \Omega',\,.\,\Omega'
\rangle$, where $\Omega'$ is in the domain of the quantum field.
This state is said to satisfy the criterion of local stability
(at $p \in M$) if there exists some monotone function $\lambda
\mapsto N_{\lambda}$ such that for any choice of finitely many
test-functions $f_1,\ldots,f_n$, the limit for $\lambda \to 0$ of
$$ \langle \Omega',\phi_{\lambda}(f_1) \cdots \phi_{\lambda}(f_n)
\Omega' \rangle  \eqno(1.2) $$
exists and equals the $n$-point function of a (non-trivial) Wightman field
over Minkowski-spacetime. (Here, the $f_j$ are identified with 
test-functions on Min\-kowski-spacetime through the chosen normal
coordinate system. Since the Wightman functions are invariant
under Lorentz-transformations, the formulated criterion is independent
of the choice of the normal coordinates.) It need, however, not be 
the case that the same state satisfies the criterion of local
stability with respect to another quantum field which generates the
same local algebras as $\phi$. Moreover, it is clear that the limiting
behaviour of $N_{\lambda}$ for $\lambda \to 0$ must be precisely known
since, for an only slightly different asymptotic behaviour of $N_{\lambda}$
near $\lambda = 0$ than the one that might lead to local stability,
the expressions (1.2) will diverge or approach 0 as $\lambda$ tends
to 0. In this case, no useful information about the short-distance
behaviour of the given theory can be gained from the scaling limits 
of (1.2).

Let us then see if we can extract some ``invariant'' meaning that may
underlie the definition of the renormalization group transformations
$R_{\lambda}$. For that purpose, we shall now specialize our
considerations to $(M,g) =$ Minkowski-spacetime. In that case, we 
have a (weakly continuous) unitary representation
$\Pg \owns (\Lx) \mapsto U(\Lx)$ of the Poincar\'e group on
$\Cal H$ which acts covariantly on the quantum field $\phi$, leaves
the (up to a phase) unique vacuum vector $\Omega$ (lying in the
domain of $\phi$) invariant, and satisfies the spectrum condition.
We write $\Aa(\Oo) := \{\phi(f) : \text{supp}(f) \subset \Oo\}''$ for
the local von Neumann algebras generated by the quantum field, and assume 
now that $(x^{\mu}) \in \Rl^4$ are coordinates of points in
Minkowski spacetime in some (arbitrary but fixed) Lorentz frame.
Then it is clear that, first of all, the renormalization group
transformations (1.2) induce mappings
$$ R_{\lambda} : \Aa(\Oo) \to \Aa(\lambda\Oo) \eqno(1.3) $$
since they act by scaling the spacetime coordinates of the quantum
field. The role played by the field strength renormalization factor
is not so immediately clear. To see what it means we note
(without going into detail here, see e.g.\ [4,6] for more information)
that it may be expected quite generally that, for 
physical models,  $N_{\lambda}$ can be suitably chosen so that the 
correlation functions (1.2) are of comparable order of magnitude
for all $\lambda$, or even converge in the scaling limit $\lambda \to 0$
and satisfy local stability. Given that this is the case,
consider the expression
\TagsOnRight
$$\split
|\langle R_{\lambda}(\phi(f))\Omega,P_{\nu}R_{\lambda}
(\phi(f))\Omega \rangle |
  &= |\langle \phi_{\lambda}(f)\Omega,[P_{\nu},\phi_{\lambda}(f)]
\Omega \rangle | \\
  &= \lambda^{-1} |\langle \phi_{\lambda}(f)\Omega,\phi_{\lambda}
(\partial_{x^{\nu}}f)\Omega \rangle |\,,
\endsplit \tag1.4$$
where $P_{\nu}$ are the generators of the translations $(\nu =
1,\ldots,4)$, i.e.\ $U(1,x) = \text{e}^{i\sum_{\nu}P_{\nu}x^{\nu}}$,
$x = (x^{\nu})$. As the correlation functions are of the same order
of magnitude for all $\lambda > 0$, we see from (1.4) that the 
energy-momentum transfer of $R_{\lambda}(\phi(f))$ (here in the
vacuum state) is approximately $\lambda^{-1}$ times the energy-momentum
transfer of $\phi(f)$. Notice that the expression (1.4) is really
a measure for the energy-momentum transfer of $R_{\lambda}(\phi(f))$,
i.e.\ the failure of $R_{\lambda}(\phi(f))$ to commute with the
$P_{\nu}$.

 Now write $\Aa := \overline{\bigcup_{\Oo}\Aa(\Oo)}^{||\,.\,||}$
and $\ax := \text{Ad}\,U(1,x)$, and define the
{\it energy-momentum transfer} $\text{EMT}(A)$ of $A \in \Aa$ as the
support of the Fourier-transform of the function
$\Rl^4 \owns x \mapsto \ax(A)$. We set $\widetilde{\Aa}(\widetilde{\Oo})
:= \{ A \in \Aa : \text{EMT}(A) \subset \widetilde{\Oo} \}$ for
$\widetilde{\Oo} \subset \Rl^4$. Then we extract from the observation
just made that the renormalization group transformations $R_{\lambda}$
induce mappings
$$ R_{\lambda} : \widetilde{\Aa}(\widetilde{\Oo}) \to
 \widetilde{\Aa}(\lambda^{-1}\widetilde{\Oo}) \eqno(1.5) $$
for each region $\widetilde{\Oo} \subset \Rl^4$ in 
``momentum space''.

It is apparent that the conditions (1.3) and (1.5) make no further
reference to a quantum field and are meaningful also in the setting
of general algebraic quantum field theory. They say essentially that
the measure of the phase-space volume which an observable $A$
occupies should not change under renormalization group transformations,
see [1,2].

 Let us recall the assumptions for a generic quantum
field theory in the operator algebraic setting [11,12], denoted
by $(\Aa,\aP,\omega)$: It is given by a local net of $C^*$-algebras
$\Oo \to \Aa(\Oo)$ indexed by the bounded open regions $\Oo$ in
Minkowski spacetime, $\Pg \owns (\Lx) \mapsto \aLx \in \text{Aut}(\Aa)$
is a covariant representation of the Poincar\'e group, and $\omega$
is the (unique) vacuum state invariant under $\aP$ and satisfying
the spectrum condition. We suppose that we are in the vacuum
representation on the vacuum Hilbertspace ${\Cal H}$ and
$\omega(\,.\,) = \langle\Omega,\,.\,\Omega\rangle$,
$\aLx = \text{Ad}\,U(\Lx)$. Moreover we emphasize that for our
approach it is quite important to assume that $\aP$ acts
{\it strongly} continuously (i.e.\ $(\Aa,\aP)$ is a $C^*$-dynamical
system). However, this is no loss of generality since one may 
always achieve this by ``smoothing'' the elements of $\Aa$ with
respect to the group action through convolution with $L^1(\Pg)$-functions,
see [1] for further discussion.

Motivated by the discussion above we now say that a family
$R_{\lambda}$, $\lambda > 0$, of maps of $\Aa$ is a renormalization
group transformation for the theory $(\Aa,\aP,\omega)$ if it
has the properties (1.3) and (1.5) and if, in addition,
$R_{\lambda}$ is continuous (in $A \in \Aa$), and uniformly bounded
in $\lambda$. This latter property expresses a certain regularity
which is reasonable in order to compare theories at different scales;
it is motivated by the analogous property of the renormalization
group transformations (1.1), but here to be interpreted in the
$C^*$-norm sense. The said conditions express the physical constraints
on renormalization group transformations in algebraic quantum field theory.
Clearly, they do not fix a particular family $R_{\lambda}$,
$\lambda > 0$, of renormalization group transformations, and
the explicit construction of renormalization group transformations is,
in general, a difficult task. However, it is one of the basic 
ingedients of algebraic quantum field theory that the interpretation
of a theory is essentially fixed by the net structure alone, so
that no renormalization group transformation fulfilling the above
conditions can be given preference over another (without introducing
additional input). Therefore, one ought to consider all of these
renormalization group transformations at an equal footing. In order
to do this, it is convenient to introduce, for a given theory
$(\Aa,\aP,\omega)$, an algebra of functions\footnote{We adopt the
convention to denote these functions, and correspondingly the objects
referring to the scaling algebra, by underlining.}
$\Rl^{+} \owns \lambda \mapsto 
\underline{A}_{\lambda}
 \in \Aa$, called {\it
scaling algebra}, which comprises all ``orbits'' $\lambda \mapsto
R_{\lambda}(A)$ of elements $A \in \Aa$ under any of the possible
renormalization group transformations.
The precise definition is as follows. 
\proclaim{1.1 Definition}
Let $(\Aa,\aP,\omega)$ be an algebraic quantum field theory
(in Minkowski spacetime). Its corresponding {\it scaling algebra}
$\underline{\Cal A}$ is the quasilocal algebra
 $\underline{\Cal A}
 := \overline{\bigcup_{\Oo}
\underline{\Cal A}
(\Oo)}^{||\,.\,||}$ of the {\it local scaling algebras}
 $\underline{\Cal A}
(\Oo)$
which are, for each bounded open region $\Oo$, defined as the sets of 
all functions $\Rl^+ \owns \lambda \mapsto
 \underline{\Cal A} \in \Aa$ with
the properties:
\par 
%% \roster
%% \item"{
$(\alpha)$ \quad
%% }"
 $||\,
\underline{A}
 \,|| := \sup_{\lambda > 0}\,||\, \underline{A}_{\lambda}\,||
  < \infty \,,$
\par
%% \item"{
$(\beta)$ \quad
%% }"
 $\underline{A}_{\lambda}
 \in \Aa(\lambda \Oo)\,,$ $\lambda > 0$,
\par 
%% \item"{
$(\gamma)$ \quad
%% }" 
$ ||\,
\underline{\alpha}_{\Lambda,x}
(\underline{A}
) - \underline{A}\,|| \to 0$ for
 $(\Lx) \to (1,0)$ \,;
%% \endroster
\par 
\noindent
here 
$$ (\underline{\alpha}_{\Lambda,x}
(\underline{A}))_{\lambda}
 : = \underline{\alpha}_{\Lambda,\lambda x}
(\underline{A}_{\lambda} )\,,
\quad \lambda > 0\,, \ \ (\Lx) \in \Pg\,. \eqno(1.6) $$
\endproclaim
Each $\underline{\Cal A}(\Oo)$ is a $C^*$-algebra with $C^*$-norm given by
$(\alpha)$, so $\Oo \to \underline{\Cal A}
(\Oo)$ is a local net of $C^*$-algebras,
and $\Pg \owns (\Lx) \mapsto 
\underline{\alpha}_{\Lambda,x} \in \text{Aut}
(\underline{\Cal A})$ is a 
strongly continuous representation of the Poincar\'e group which
acts covariantly on the net of scaling algebras.

To comment on the meaning of the points $(\alpha)$-$(\gamma)$ in the
definition of the scaling algebra, suppose that $R_{\lambda}$,
$\lambda > 0$, is a renormalization group transformation
for $(\Aa,\aP,\omega)$, and set
 $\underline{A}_{\lambda}
 := R_{\lambda}(A)$ for some
$A \in \Aa(\Oo)$. Then $(\beta)$ is just the condition (1.3),
while $(\alpha)$ expresses the uniform boundedness of $R_{\lambda}$
in $\lambda$. It is not so obvious that (1.5) amounts to demanding
$(\gamma)$, but this has been shown in [1, Lemma 3.1] (for the
action of the translations; the analogous requirement for the 
Lorentz-transformations restricts the behaviour of the 
angular-momentum transfer under renormalization group 
transformations, cf.\ [1].)

In summary, the scaling algebra is to be viewed as formed
by the collection of all orbits $\lambda \mapsto R_{\lambda}(A)$
of elements $A \in \Aa$ under any ``abstract'' renormalization
group transformation characterized by the conditions given
above, particularly (1.3) and (1.5). The major advantage
of looking at renormalization group transformations in that way
is of course that we can use the powerful operator algebraic
machinery to analyze them and to consequently obtain results, as
will become clear in the next section.
\heading
2. Scaling Limits
\endheading
Let $(\Aa,\aP,\omega)$ be an algebraic quantum field theory on
Minkowski spacetime and 
$\underline{\Cal A}$
 its scaling algebra together
with the lifted action
 $\underline{\alpha}_{{\Cal P}^{\uparrow}_+}$,
 cf.\ (1.6). We will in the
following speak of states on $\Aa$ which are locally normal to
the vacuum $\omega$ as {\it physical states} of the given
theory $(\Aa,\aP,\omega)$. With each such physical state
$\omega'$ one may associate a family 
$(\underline{\omega}_{\lambda}'
 )_{\lambda > 0}$
of states on $\underline{\Cal A}$, the ``scaled lifts'' of $\omega'$,
which are defined by 
$$ \underline{\omega}_{\lambda}'
(\underline{A}
) := \omega'
(\underline{A}_{\lambda})\,, \quad \lambda > 0\,,\ \ 
\underline{A} \in 
\underline{\Cal A}
\,.\eqno(1.7) $$
We mention that, if $(\underline{\Cal H}
,\underline{\pi}_{\lambda}
,\underline{\Omega}_{\lambda}
)$ is the GNS-representation
of $\underline{\omega}_{\lambda}$
, and
 ${}^{(\lambda)}\!\underline{\alpha}_{{\Cal P}^{\uparrow}_+}$
 denotes the induced action of the
Poincar\'e group in the GNS-representation, then the theory
$(\underline{\pi}_{\lambda}
(\underline{\Cal A}
),
{}^{(\lambda)}\!\underline{\alpha}_{{\Cal P}^{\uparrow}_+}
,
\underline{\omega}_{\lambda}
)$ is isomorphic to the ``scaled''
version $(\Aa_{\lambda},\alpha^{(\lambda)}_{\Pg},\omega)$ of the
given theory, where $\Aa_{\lambda}(\Oo) := \Aa(\lambda \Oo)$
and $\alpha^{(\lambda)}_{\Lx} := \alpha_{\Lambda,\lambda x}$.
(See [1]. A similar result holds for physical states.)
Here we recall that two theories $(\Aa^{(1)},\alpha^{(1)}_{\Pg},
\omega^{(1)})$ and $(\Aa^{(2)},\alpha^{(2)}_{\Pg},\omega^{(2)})$
 are called isomorphic if there is a $C^*$-algebraic isomorphism
$\rho : \Aa^{(1)} \to \Aa^{(2)}$ which (1) preserves the
net-structure, $\rho(\Aa^{(1)}(\Oo)) = \Aa^{(2)}(\Oo)$ for all
regions $\Oo$, (2) intertwines the Poincar\'e group actions,
$\rho \circ \alpha^{(1)}_{\Lx} = \alpha^{(2)}_{\Lx} \circ \rho$
for all $(\Lx) \in \Pg$, and (3) connects the vacua,
$\omega^{(1)} = \omega^{(2)} \circ \rho$. It is clear that 
isomorphic theories are physically indistinguishable, they
describe identical physical situations.
Thus the states (1.7) on the scaling algebra carry for each $\lambda$
the same information as the physical states on the by the factor
$\lambda$ ``scaled'' version of the original theory.
Therefore, when one formally proceeds to the limit $\lambda \to 0$
in $(1.7)$ then the extreme short distance remnants of the
originally given theory, so to speak the processes that take place
at ``zero spatio-temporal scale'', appear.

There occurs now a difficulty since the limits of (1.7) for
$\lambda \to 0$ need in general not exist. However, when we
view the family 
$(\underline{\omega}_{\lambda}'
)_{\lambda > 0}$ of states on 
$\underline{\Cal A}$
as a net of states indexed by the positive reals directed towards 0,
then the Banach-Alaoglu theorem asserts the existence of 
weak-* limit points of this net.
\proclaim{2.1 Definition}
The set of weak-* limit points of 
$(\underline{\omega}_{\lambda}'
)_{\lambda > 0}$ for
$\lambda \to 0$ is denoted by $SL(\omega') = 
\{\underline{\omega}'{}_{0,\iota}
 : \iota \in
{\Bbb I} \}$ (where ${\Bbb I}$ is some abstract index set 
labelling the weak-* limit points). Each 
$\underline{\omega}'{}_{0,\iota}$,
 $\iota \in
{\Bbb I}$, is called a {\it scaling limit state} of $\omega'$.
\endproclaim
Drawing on the fact that $\bigcap_{\Oo \owns 0} \Aa(\Oo)^- =
{\Bbb C}1$ (see [13], and also [1]), one obtains that
$$ \lim_{\lambda \to 0} \, |
\underline{\omega}_{\lambda}'(
\underline{A}_{\lambda}
) - 
\underline{\omega}_{\lambda}
(\underline{A}_{\lambda}) | = 0\,,
\quad 
\underline{A} \in 
\underline{\Cal A}  \,,$$
for any physical state $\omega'$ on $\Aa$ (see [1] for
the proof, which is based on an argument by Roberts [14]).
Hence the scaling limit states of any physical state
coincide with the scaling limit states of the vacuum
 and so it suffices in the following to consider only
the latter.
We denote by $(\underline{\Cal H}_{0,\iota}
,\underline{\pi}_{0,\iota},
\underline{\Omega}_{0,\iota}
)$ the GNS-representation of
$\underline{\omega}_{0,\iota}
 \in SL(\omega)$. Obviously
 $\underline{\omega}_{0,\iota}$
 is invariant under the
lifted action of the Poincar\'e group, and we can pass to
the next
\proclaim{2.2 Definition}
We write ${\Cal A}_{0,\iota}(\Oo)
 := \underline{\pi}_{0,\iota}(\underline{\Cal A}
(\Oo))$, and $\alpha^{(0,\iota)}_{{\Cal P}^{\uparrow}_+}$
 for the
induced action of the Poincar\'e group in the GNS-representation
of $\underline{\omega}_{0,\iota}$
, so $\alpha^{(0,\iota)}_{\Lambda,x}
 \circ \underline{\pi}_{0,\iota}
 = \underline{\pi}_{0,\iota}
 \circ 
\underline{\alpha}_{\Lambda,x}$.
 We call
$\Oo \to {\Cal A}_{0,\iota}
(\Oo)$ the {\it scaling limit net} and
$({\Cal A}_{0,\iota}
,\alpha^{(0,\iota)}_{{\Cal P}^{\uparrow}_+},
\underline{\omega}_{0,\iota})$ the {\it scaling limit theory} corresponding to
the scaling limit state 
$\underline{\omega}_{0,\iota}
 \in SL(\omega)$.
\endproclaim
We quote some fairly immediate results about the scaling
limit theories which are associated with each scaling limit
state $\underline{\omega}_{0,\iota}
 \in SL(\omega)$ from [1].
\proclaim{2.3 Proposition}
(a) Each scaling limit net $\Oo \to
 {\Cal A}_{0,\iota}
(\Oo)$ fulfills the
locality condition. 
\par 
\noindent
(b) Each scaling limit theory $(
{\Cal A}_{0,\iota},
\alpha^{(0,\iota)}_{{\Cal P}^{\uparrow}_+},
\underline{\omega}_{0,\iota}
)$ is a theory
over Minkowski spacetime in vacuum representation:
$\alpha^{(0,\iota)}_{{\Cal P}^{\uparrow}_+}$
 is strongly continuous, acts covariantly on the scaling limit
net, leaves 
$\underline{\omega}_{0,\iota}$
 invariant, and the spectrum condition is
satisfied so that 
$\underline{\omega}_{0,\iota}$
 is a vacuum state. 
\par 
\noindent
(c) $\underline{\omega}_{0,\iota}$
 is a pure state (this holds for spacetime-dimension $\ge 2$).
 \par 
\noindent
(d) Suppose that geometric modular action for wedge regions
holds in the 
\linebreak
 originally given theory, i.e.\ with $\Delta,J$ the
modular objects of $\Aa(W)^-,\Omega$, and
$W := \{ x \in \Rl^4: 0 < |x^0| < x^1 \}$ a wedge region, we 
posit that 
$$ J\Aa(\Oo)^-J = \Aa(j(\Oo))^- \quad \text{and} \quad
\text{\rm Ad}\,\Delta^{it} = \alpha_{\Lambda_{2\pi t}}\,,\ \ t \in \Rl\,,
\eqno(1.8) $$
where $j(x^0,x^1,x^2,x^3) = (-x^0,-x^1,x^2,x^3)$ and the Lorentz
boost $\Lambda_s$ acts like
 $\left( {}^{{}\,\,\text{cosh}(s)\ \,\,\,-\text{sinh}(s)}
_{-\text{sinh}(s)\ \ \,\text{cosh}(s)} \right)$ on the first two
coordinates, leaving the others fixed. 
\par
Then the relations $(1.8)$ hold in each scaling limit theory,
with $\Aa$ replaced by $\Aoi$ and $\Delta_{0,\iota},J_{0,\iota}$
defined as modular objects associated with 
${\Cal A}_{0,\iota}
(W)^-,
\underline{\Omega}_{0,\iota} $.
\endproclaim
Let us just recapitulate: We have associated with the originally
given theory $(\Aa,\aP,\omega)$ --- which we will from now on
often refer to as the {\it underlying theory} ---  the scaling algebra,
and the families of scaled lifts of physical states, whose
weak-* limit points for $\lambda \to 0$ we have collected in the
set of scaling limit states $SL(\omega)$. Then we passed
to the scaling limit theories and have thus assigned to the 
underlying theory a whole family of scaling limit theories
$(
{\Cal A}_{0,\iota}
,
\alpha^{(0,\iota)}_{{\Cal P}^{\uparrow}_+}
,
\underline{\omega}_{0,\iota})$, $\iota \in {\Bbb I}$, which encodes
information about the short-distance properties of the underlying
theory.

To see what information we gain in that way, a natural first step
is to identify the members in the family of scaling limit theories
which are isomorphic, and to form the corresponding isomorphy
classes. Then the following mutually exclusive cases can occur.

\roster
\item"({\bf C})" There is only one isomorphy class of scaling limit
theories, and  
\newline
${\Cal A}_{0,\iota} = {\Bbb C}1.$ We call this case the
{\it classical} (or {\it trivial}) scaling limit.
\item"({\bf Q})" There is only one isomorphy class of scaling limit
theories, and 
${\Cal A}_{0,\iota} 
$ is non-Abelian. This case will be called
the {\it quantum} scaling limit. 
\item"({\bf D})" There is more than one isomorphy class of scaling
limit theories. We refer to this case as {\it degenerate} scaling
limit.
\endroster

If the underlying theory has a classical scaling limit, this
means that its phase-space behaviour is rapidly worsening
at small scales in the sense that, as $\lambda \to 0$,
the expectation values of non-Abelian elements in the algebras
$\Aa(\lambda\Oo)$ whose transferred energy-momentum scales like
$\lambda^{-1}$ vanish for physical states.
On the other hand, if the underlying theory has a quantum scaling
limit, this is a sign of a stable behaviour of its dynamics
at small scales and corresponds to the situation of local
stability. Physical theories, like QCD, are believed to
belong to this class. In contrast to that, an underlying theory
with a degenerate scaling limit has a very irregular and unstable
dynamical behaviour at small scales. It is expected that this
case occurs for theories which in the terminology of perturbative
quantum field theory do not posses an ultra-violet fixed point
under renormalization group transformations.

It should be pointed out that for each of the alternatives
({\bf C}),({\bf Q}),({\bf D}) there are underlying theories having
the respective behaviour in the scaling limit. Examples for
theories satisfying ({\bf C}) or ({\bf D}) are constructed
from generalized free scalar fields [15]. In spacetime dimension = 3,4
the free scalar field with mass $\ge 0$ is an example of an
underlying theory with quantum scaling limit, with all scaling limit
theories being isomorphic to the free, massless scalar field [16].

There is a general result connecting the existence of dilations
as geometrical symmetries and the non-degeneracy of the scaling
limit, which will be quoted next from [1]. We say that
a theory $(\Aa,\aP,\omega)$ is dilation covariant with dilation
invariant vacuum if there is a family $\delta_{\mu} \in \text{Aut}(\Aa)$,
$\mu > 0$, with $\delta_{\mu}(\Aa(\Oo)) = \Aa(\mu \Oo)$,
$\delta_{\mu} \circ \aLx = \alpha_{\Lambda,\mu x} \circ \delta_{\mu}$,
and $\omega \circ \delta_{\mu} = \omega$.

\proclaim{2.4 Proposition} 
(a) Suppose that the underlying theory is dilation covariant with
dilation invariant vacuum and satisfies, moreover, the Haag-Swieca
compactness condition $[17,18]$. Then the scaling limit of this
theory is non-degenerate, and all scaling limit theories are
isomorphic to the underlying theory itself. 
\par 
\noindent
(b) Assume that the scaling limit of an underlying theory
is nondegenerate (case $(${\bf C}$)$ or $(${\bf Q}$)$). Then each of the
(isomorphic) scaling limit theories is dilation covariant with
dilation invariant vacuum.
\endproclaim
Notice that (a) can be read as saying that dilation covariant theories
with dilation invariant vacuum are ``fixed points'' under
renormalization group transformations.
It is worth mentioning that conditions on the nature of the
scaling limit, particularly a distinction of the cases ({\bf C})
and ({\bf Q}), can be formulated in terms of the energy-nuclearity or
-compactness behaviour of the underlying theory. The reader
is referred to [3] for a survey of these matters.

Since one of the motivations for studying renormalization group
transformations within the abstract approach to local quantum field
theory was the incomplete conceptual basis of understanding of
confinement, let us now outline, following [4], how this issue
can be addressed in the setting presented so far. To this 
end we must recollect a few basic notions and results of the
Doplicher-Haag-Roberts framework for the description of charges
via super\-selection sectors in algebraic quantum field theory
[11,19]. One starts with a theory $(\Aa,\aP,\omega)$
satisfying Haag-duality, and imposes on the set of physical
states a further selection criterion expressing the physical
properties of the charged states of the given theory, typically
by demanding that they be normal to the vacuum outside of 
bounded regions [19] or outside of spacelike cones [20].
The equivalence classes of all pure such states are called
superselection sectors, and the set of all of them, the charge
spectrum, will be denoted by $\Sigma$. Doplicher and Roberts [21]
were able to prove that then there is a Hilbertspace
$\Hh^{\Ff}$ containing the vacuum Hilbertspace $\Hh = (\Aa\Omega)^-$,
on which a $C^*$-algebraic net of field algebras $\Oo \to
\Ff(\Oo) \supset \Aa(\Oo)$ acts, together with an extension
of the covariant automorphic action of $\Pg$ to the net of
field algebras. In addition, there is a compact group $G$ (the
global gauge group), and a representation $G \owns g \mapsto V(g)$
by unitaries on $\Hh^{\Ff}$, so that $\text{Ad}\,V(g)$ acts as
an internal symmetry on the net of field algebras,
$V(g)\Omega = 0$, and $\Aa(\Oo)$ consists precisely of the
fixed points in $\Ff(\Oo)$ under the action of $\text{Ad}\,V(g)$.
The spectrum of $G$ coincides with $\Sigma$, and $\Hh^{\Ff} =
(\Ff\Omega)^-$.

A similar construction can be carried out for the scaling
limit theory. (With the provision that for this purpose the
superselection theory may have to be extended to comprise
massless particles with non-local charges.) Let us assume that 
the underlying theory has a quantum scaling limit. Then, since
all scaling limit theories are isomorphic, we can denote the
scaling limit theory just by $({\Cal A}_0,\alpha^{(0)}_{\Pg},\omega_0)$,
dropping the labelling index $\iota$. We denote by
$\Ff^0$, $G_0$ and $\Sigma_0$ the field algebra, gauge group,
and charge spectrum, respectively, of the scaling limit theory.
The particles and charges described by these entities are
called {\it ultraparticles} and {\it ultracharges} [2,4].
We write $\Ff_0$ for the scaling limit algebra of the field algebra
of the underlying theory --- defined by first constructing
a scaling field algebra $\underline{\Ff}$ in analogy to the
scaling algebra of the observables $\underline{\Aa}$, defining scaled
lifts of locally normal states on $\Ff$, and then passing to the
GNS-representation of $\underline{\Ff}$ corresponding to the 
scaling limits states. It is obvious that
$\Aa_0(\Oo) \subset \Ff_0(\Oo)$ and $\Aa_0(\Oo) \subset \Ff^0(\Oo)$,
but in general, $\Ff_0(\Oo)$ and $\Ff^0(\Oo)$ need not 
coincide. One may expect that $\Ff_0(\Oo) \subset \Ff^0(\Oo)$,
and it may happen that $\Hh^{\Ff_0} = (\Ff_0\Omega)^-$ is a proper
subspace of $\Hh^{\Ff^0} = (\Ff^0\Omega)^-$. According to
[4], this can be taken as criterion for {\it confinement}.
Any vector state $\langle \psi,\,.\,\psi \rangle$ with 
$\psi \in \Hh^{\Ff^0} \ominus \Hh^{\Ff_0}$ is identified with
the state of an ultraparticle which is confined, since that
state is not created by the field algebra of the underlying theory
at any finite scale.

 It is worth pointing out that here one has
reached at a criterion for confinement is entirely independent
of any ``coordinatization'' of the underlying theory
(or, for that matter, the scaling limit theory) by quantum
fields, and so the occurrence or non-occurrence of confinement
is an invariant of the given theory, completely based on a formulation
in terms of observable quantities. The reader is referred to
[4] for further discussion. There exists an example for a 
theory where confinement in the just described sense occurs:
It is given by the Schwinger model of massless quantum electrodynamics
in two-dimensional Minkowski spacetime. The net of local observables
of this model coincides with that of the free, neutral massive scalar field.
It is known that there are no charged sectors which are locally
normal to the vacuum for this model [22]. 
But it can be shown that in the scaling limit
theories there appear charged sectors, i.e. superselection sectors
which are locally, but not globally normal to the scaling limit 
vacua [4,16]. This fact indicates that the scaling algebra approach
appears to be quite promising in having the potential to lead
to a satisfactory conceptual understanding of the confinement problem.  
\heading
3. Local Stability
\endheading
It was mentioned in Section 1 that in quantum field theory on curved
\linebreak
spacetimes there is no straightforward specification of physical
states as those locally normal to a vacuum state since for
the latter there will in general, due to the absence of spacetime
symmetries, be no candidate. Following [9,10,11], one
can impose the principle of local stability as a constraint on
physical states. In order to be able to give an intrinsic definition
of local stability in the spirit of the discussion in the preceding
sections, one is naturally inclined to use the scaling algebra framework
also for algebraic quantum field theories in curved spacetime.
Again, the absence of symmetries entails some complications since
it means that there is no immediate generalization of the condition
$(\gamma)$ constraining the phase-space behaviour under 
renormalization group transformations. However, under quite general
conditions one can give appropriate versions of $(\gamma)$ also
for theories in curved spacetimes. We shall fix these conditions now,
and assume that $(M,g)$ is a spacetime which is globally hyperbolic
([23,24], we make this assumption mainly for convenience here,
it is not essential), and that ${\Cal R} \to \Aa({\Cal R})$ is a 
local and primitively causal net of $C^*$-algebras indexed by
the relatively compact open subsets of $M$; we again write
$\Aa = \overline{\bigcup_{\Cal R}\Aa({\Cal R})}^{||\,.\,||}$.
Locality means that the algebras $\Aa({\Cal R}_1)$ and $\Aa({\Cal R}_2)$
commute elementwise if ${\Cal R}_1 \subset {\Cal R}_2^{\perp}$ where
the causal complement $S^{\perp}$ of $S \subset M$ is defined as the
largest open set of  points in $M$ which cannot be connected to
$S$ by any causal curve. Primitive causality, a strong form of the
``time-slice axiom'', means that if $G$ is an acausal hypersurface
in $M$ (no pair of different points in $G$ can be joined by a causal
curve), then it holds that 
$$ \bigcap_{{\Cal R} \supset G} \Aa({\Cal R}) =
\Aa({\Cal R}_G) \quad \text{with}\quad {\Cal R}_G := 
(G^{\perp})^{\perp} \,. $$
In view of primitive causality one can, whenever one is given
a foliation $F : \Rl \times \Sigma \to M$ in Cauchy-surfaces
(i.e.\ $\Sigma$ is a 3-dim.\ smooth manifold, $F$ a diffeomorphism,
and for each $t \in \Rl$, $F(\{t\} \times \Sigma)$ is a Cauchy-surface;
$t$ has the significance of a time parameter), introduce canonically
the net of algebras ``at foliation time $t$'',
$$ \Sigma \supset G \mapsto \Aa_{F(t,G)} := \Aa({\Cal R}_{F(\{t\} \times
G)}) \,. $$
Then it makes sense to assume that there is a dynamical time-evolution
for this foliation, given by a {\it propagator family}
$\alpha^{(F)}_{t,t'} \in \text{Aut}(\Aa)$, $t \ge t' \in \Rl$, satisfying
the propagator identity
$$ \alpha^{(F)}_{t,t'} \circ \alpha^{(F)}_{t',s} =
\alpha^{(F)}_{t,s}\,, \quad \alpha^{(F)}_{t,t} = \text{id}_{\Aa} \,, $$
and preserving the ``localization at equal time'',
$$ \alpha^{(F)}_{t,t'}({\Cal A}_{F(t',G)}) = \Aa_{F(t,G)}\,,
\quad G \subset \Sigma\,. $$
We note that the $C^*$-algebraic nets of local observables for
 fields obeying a linear hyperbolic equation of motion in a globally
hyperbolic spacetime typically are primitively causal and admit 
propagator families for arbitrary foliations in Cauchy-surfaces [5,25,26].

Some foliations have distinguished poperties.   
For instance, given $p \in M$, it can be shown that there are
foliations which near $p$ are geodesic (i.e.\ the curves
$t \mapsto F(t,q)$ are geodesics) and hypersurface orthonormal
(i.e.\ $\frac{\partial}{\partial t} F(t,q)$ equals the
(future- or past-pointing) unit-normal vector of $F(\{t\} \times \Sigma)$
at $F(t,q)$). Moreover, one can always find such foliations
with the additional property that the normal at $p$ is parallel
to any prescribed timelike vector $v \in T_pM$. The significance
of such geodesic, hypersurface orthonormal foliations is that for them
$\frac{\partial}{\partial t}F(t,q)$ is the curved spacetime analog
of a non-rotating
linear velocity field  [24].

The basic assumption we will therefore make is that, for the underlying
theories which we consider, there is around any point $p \in M$
a collection $(\alpha^{(F)})_{F \in \Phi}$ of propagator families,
where the set $\Phi$ contains geodesic, orthonormal foliations
around $p$ all of whose normals at $p$ fill the cone of timelike
vectors in $T_pM$. The point $p$ is supposed to lie on the
foliation-time $= 0$ surfaces.
 We shall furthermore suppose that each
$\alpha^{(F)}_{t,t'}$ is strongly continuous in $t,t'$. We
call $(\alpha^{(F)})_{F \in \Phi}$ a dynamics (at $p$). In 
general there may be several choices for such a dynamics. We take
here the point of view that a particular one has been selected, as
part of the specification of the theory.

Now we can define the scaling algebra at a given $p \in M$.
Pick a normal coordinate chart $u = (x^{\mu})$ at $p$ (recall that
$u(p) =0$). The coordinate basis of that chart allows to
canonically identify $T_pM$ with Minkowski spacetime, and simultaneously,
under this identification, $u$ takes values in $T_pM$.

 Then we define
the {\it scaling algebra} at $p$ as the quasilocal algebra
$\underline{\Aa} = \overline{\bigcup_{\Oo} \underline{\Aa}(\Oo)}^{||\,.\,||}$
of the local scaling algebras $\underline{\Aa}(\Oo)$, which for 
$\Oo \subset T_pM$ open and bounded consist of the functions
$\Rl^+ \owns \lambda \mapsto \underline{A}_{\lambda} \in \Aa$,
subject to the following requirements:
\roster
\item"$(\alpha')$" $||\,\underline{A}\,|| :=
\sup_{\lambda > 0}\,||\,\underline{A}_{\lambda}\,|| < \infty$ ,
\item"$(\beta')$" $\underline{A}_{\lambda} \in \Aa(u^{-1}(\lambda \Oo
\cap \text{Ran}(u)))$, \quad $\lambda > 0$,
\item"$(\gamma')$" $\limsup_{\lambda\to 0}\,
||\,\alpha^{(F)}_{\lambda t,\lambda t'}(\underline{A}_{\lambda})
- \underline{A}\,|| \to 0$\ \ \  for\ \ \  $t,t' \to 0$ \newline
and $\underline{\alpha}^{(F)}_{t,t'} \in \text{End}(\underline{\Aa})$,
%%%%%%% \newline
for all propagators $\alpha^{(F)}_{t,t'}$ of the dynamics at $p$,
\endroster
where
$$ (\underline{\alpha}^{(F)}_{t,t'}(\underline{A}))_{\lambda}
:= \alpha^{(F)}_{\lambda t,\lambda t'}(\underline{A}_{\lambda}) \,. $$
\par
Some comments are in order here. For $(\beta')$ note that, since
$u(p) = 0$, $\lambda \Oo$ will be in the range of $u$ for 
sufficiently small $\lambda$. It is understood that
$\underline{A}_{\lambda}$ takes values in ${\Bbb C} 1$ if
$\lambda\Oo \cap \text{Ran}(u) = \emptyset$. The
$\underline{A}_{\lambda}$ become, for $\lambda \to 0$, localized
at $p$; so $(\beta')$ plays the role of $(\beta)$ in Definition 1.1,
here with respect to the point $p \in M$. Likewise, $(\gamma')$
is our curved spacetime version of $(\gamma)$ in Definition 1.1,
restricting the behaviour of the energy-momentum transfer under
renormalization group transformations.
Besides depending on the chosen point $p \in M$, the definition
of the scaling algebra depends also on the choice of the
normal coordinate system $u = (x^{\mu})$. But this dependence becomes
trivial in the scaling limit, owing to the fact that for any pair
$u,u'$ of normal coordinate charts at $p$ it holds that
$u'\circ u^{-1}(x) = \Lambda x + O(x^2)$ near $x = 0$ with
some Lorentz-transformation $\Lambda$.

Now let a state $\omega$ on $\Aa$ be given. We shall use the notation
$\Aa_{\omega}({\Cal R}) := \pi_{\omega}(\Aa({\Cal R}))^-$ for the local
von Neumann algebras in the GNS-representation of $\omega$, and we 
will only consider those states for which
$$ {\Cal H}_{\omega}\ \ \  \text{is separable,\ \  and}\ \ \ 
\bigcap_{{\Cal R} \owns p}\Aa_{\omega}({\Cal R}) = {\Bbb C}1 \,.
\eqno(3.1) $$
In the same manner as in Section 2, we obtain the scaled lifts of
$\omega$,
$$ \underline{\omega}_{\lambda}(\underline{A}) = \omega(\underline{A}
_{\lambda})\,, \quad \underline{A} \in \underline{\Aa}\,, $$
and the corresponding weak-* limit points $\underline{\omega}_{0,\iota}
\in SL(\omega)$ for $\lambda \to 0$, called scaling limit states
of the scaling algebra at $p$.

For the next proposition it is convenient to focus on the
inner regularized scaling limit nets
$$ \Oo \to \Aa_{0,\iota}(\Oo) := \overline{
\bigcup_{\overline{\Oo_1} \subset \Oo}
 \underline{\pi}_{0,\iota}(\underline{\Aa}(\Oo_1)) }^{||\,.\,||}\,,$$
indexed by the bounded open regions $\Oo \subset T_pM$, where
$\underline{\pi}_{0,\iota}$ is the GNS-representation of
$\underline{\omega}_{0,\iota} \in SL(\omega)$. Again we write
$\Aa_{0,\iota} := \overline{\bigcup_{\Oo} \Aa_{0,\iota}(\Oo)}^{||\,.\,||}$.
\proclaim{3.1 Proposition}
(a) $\Oo \to \Aa_{0,\iota}(\Oo)$ is a local and primitively causal
$C^*$-algebraic net on Minkowski spacetime for all scaling limit states
$\underline{\omega}_{0,\iota}$ of the scaling algebra at $p$.
\par
\noindent
(b) If $\text{\rm ker}\,\underline{\pi}_{0,\iota}$ is invariant under the
action of $\underline{\alpha}^{(F)}_{t,t'}$, then there is a
strongly continuous 
propagator family $\alpha^{(F;0,\iota)}_{t,t'} \in \text{\rm Aut}
(\Aa_{0,\iota})$, $t \ge t'$, so that
$$ \alpha^{(F;0,\iota)}_{t,t'} \circ \underline{\pi}_{0,\iota}
= \underline{\pi}_{0,\iota} \circ \underline{\alpha}^{(F)}_{t,t'}\,.$$
Moreover, when $v$ is the normal vector of the foliation $F$ at $p$,
then this propagator family acts as ``equal-time translation in 
$v$-direction'', i.e.\ we have
$$ \alpha^{(F;0,\iota)}_{t,t'} (\Aa_{0,\iota}(\Oo + t'v)) =
 \Aa_{0,\iota}(\Oo + tv) $$
for all double cones $\Oo$ in Minkowski spacetime based
on the hyperplane passing through $0$ which is orthogonal to $v$.
\endproclaim
The just stated result (see [5] for a proof) motivates the 
following version of local stability for a theory on a curved
spacetime.
\proclaim{3.2 Definition} Let $\omega$ be a state on $\Aa$ obeying (3.1).
We say that a scaling limit state $\underline{\omega}_{0,\iota} \in
SL(\omega)$ on the scaling algebra at $p$ fulfills the
{\it condition of local stability} if: $\underline{\pi}_{0,\iota}$
is irreducible, $\text{\rm ker}\,\underline{\pi}_{0,\iota}$ is invariant
under $\underline{\alpha}^{(F)}_{t,t'}$ for all $\alpha^{(F)}$ of the
dynamics at $p$, and there exists a strongly continuous
covariant representation 
$\Rl^4 \owns x \mapsto \alpha^{(0,\iota)}_x \in \text{\rm Aut}(\Aa_{0,\iota})$
of the translations which leaves $\underline{\omega}_{0,\iota}$
invariant, satisfies the spectrum condition, and with the property
$$ \alpha^{(F;0,\iota)}_{t + \tau,\tau} = \alpha^{(0,\iota)}_{tv}\,,
\quad t \ge 0,\ \tau \in \Rl\,, \eqno(3.2) $$
where $v$ is the normal vector of $F$ at $p$.
\endproclaim
In other words, local stability of $\underline{\omega}_{0,\iota}$
means that $(\Aa_{0,\iota},\alpha^{(0,\iota)}_x,\underline{\omega}_{0,\iota}
)$ gives, apart from Lorentz-transformations, a theory on
Minkowski spacetime, where the action of the translations is 
asymptotically related to the dynamics of the underlying theory
through (3.2).

We finally quote from [5] the curved spacetime generalization of a
result proved in [1] connecting the type III$_1$-property of the
von Neumann algebras of the underlying theory with the existence 
of a state fulfilling local stability for a non-trivial scaling limit
in which also wedge-duality is realized. It extends the arguments 
by Fredenhagen [8] for the Wightman-field setting  (see also
[27] for the corresponding curved spacetime version) to our
scaling algebra approach.
\proclaim{3.3 Proposition}
Let ${\Cal R} = (G^{\perp})^{\perp}$, where $G$ is an open subset
of some Cauchy-surface so that $\partial G$ is piecewise smooth and
$G^{\perp} \ne \emptyset$. Let $\omega$ be a state on $\Aa$ 
satisfying $(3.1)$, and suppose that there are a point $p$ in the 
smooth part of $\partial G$ and a scaling limit state $\underline{\omega}
_{0,\iota}$ on the scaling algebra at $p$ fulfilling local stability,
non-triviality, i.e.\ $\Aa_{0,\iota} \ne {\Bbb C}1$, and 
wedge-duality, i.e.\
$$ \Aa_{0,\iota}(W)' = \Aa_{0,\iota}(W^{\perp})'' $$
for all wegde-regions $W$ in Minkowski spacetime, and with the causal
complement taken here with respect to the Minkowski metric.

Then the von Neumann algebra $\Aa_{\omega}({\Cal R})$ of the underlying
theory is of type III$_1$.
\endproclaim
\heading
Concluding Remarks
\endheading
The scaling algebra approach to renormalization group transformations
provides the means for a 
model-independent, intrinsic analysis of the extreme short-distance behaviour
 of a given
algebraic quantum field theory. We have pointed out that there are
promising indications that this method will lead to a satisfactory
conceptual understanding of the notion of confined charges and
particles used in elementary particle physics in terms of
ultracharges and ultraparticles. One may also hope that the framework
sheds some light on the question of how 
the concept of local gauge transformations
can be understood within the setting of algebraic quantum field theory.

For the situation of algebraic quantum field theory in curved
spacetime we have indicated that the scaling algebra approach
gives a model-independent notion of local stability, which may
be taken as a selection criterion for physical states. The
setting presented here is  in several respects incomplete,
and we hope to address that matter more profoundly elsewhere.
\heading
Acknowledgements
\endheading
I would like to thank D.\ Buchholz for the opportunity of joint
work on the topic presented here, as well as for numerous discussions
about it. I also owe thanks to the members of the operator algebras
group at the Dipartimento di Matematica at Tor Vergata for their
kind hospitiality during my stay in Rome, which was funded by a
Von Neumann Fellowship of the Operator Algebras Network,
EC Human Capital and Mobility Programme.
\heading
References
\endheading

\item{[1]} D.~Buchholz, R.~Verch: {\it Scaling algebras and renormalization
group in algebraic quantum field theory}, Rev.~Math.~Phys.\ {\bf 7}
(1995) 1195
\item{[2]} D.~Buchholz: {\it On the manifestations of particles}. In:
R.~N.~Sen, A.~Gersten (eds.): {\it Mathematical Physics Towards the
21st Century}, Ben Gurion University Press, Beer-Sheva (1994)
\item{[3]} D.~Buchholz: {\it Phase-space properties of local observables
and structure of scaling limits}, Ann.~Inst.~H.~Poincar\'e {\bf 64}
(1996) 433
\item{[4]} D.~Buchholz: {\it Quarks, gluons, colour: Facts or fiction?},
 Nucl.~Phys.\ {\bf B 469} (1996) 333
\item{[5]} R.~Verch, PhD.~Thesis, University of Hamburg, 1996; and
work in preparation
\item{[6]} P.~Becher, M.~B{\"o}hm, H.~Joos: {\it Gauge Theories of
Strong and Electroweak Interactions}, Wiley, New York (1984);
\newline
M.~LeBellac: {\it Quantum and Statistical Field Theory}, Oxford 
University Press, Oxford (1991);
\newline
J.~Zinn-Justin: {\it Quantum Field Theory and Critical Phenomena},
Clarendon Press, Oxford (1989)
\item{[7]} H.~J.~Borchers: {\it {\"U}ber die Mannigfaltigkeit der
interpolierenden Felder zu einer kausalen S-matrix}, Nuouvo
Cimento {\bf 15} (1960) 784
\item{[8]} K.~Fredenhagen: {\it On the modular structure of local observables},
Commun.
\linebreak
Math.~Phys.\ {\bf 97} (1985) 79
\item{[9]} R.~Haag, H.~Narnhofer, U.~Stein: {\it On quantum field theory
in gravitational background}, Commun.~Math.~Phys.\ {\bf 94} (1984) 219
\item{[10]} K.~Fredenhagen, R.~Haag: {\it Generally covariant quantum
field theories and scaling limits}, Commun.~Math.~Phys.\ {\bf 108}
(1987) 91
\item{[11]} R.~Haag: {\it Local Quantum Physics}, 2nd  edn.,
Springer Verlag, Berlin-Heidel\-berg-New York (1996)
\item{[12]} R.~Haag, D.~Kastler: {\it An algebraic approach to
quantum field theory}, J.~Math.
\linebreak
Phys.\ {\bf 4} (1964) 848
\item{[13]} A.~S.~Wightman: {\it La th\'eorie quantique locale
et la th\'eorie quantique des champs}, Ann.~Inst.~H.~Poincar\'e
{\bf 1} (1964) 403
\item{[14]} J.~E.~Roberts: {\it Some applications of dilatation
invariance to structural questions in the theory of local
observables}, Commun.~Math.~Phys.\ {\bf 37} (1974) 273
\item{[15]} D.~Buchholz, M.~Lutz, S.~Mohrdieck: In preparation
\item{[16]} D.~Buchholz, R.~Verch: In preparation
\item{[17]} R.~Haag, J.~A.~Swieca: {\it When does a quantum field thoery
describe particles?}, Commun.~Math.~Phys. {\bf 1} (1965) 308
\item{[18]} D.~Buchholz, M.~Porrmann: {\it How small is the phase space
in quantum field theory?}, Ann.~Inst.~H.~Poincar\'e {\bf 52} (1990)
631
\item{[19]} S.~Doplicher, R.~Haag, J.~E.~Roberts: {\it Fields,
observables and gauge transformations}, Commun.~Math.~Phys. {\bf 13}
(1969) 1; Commun.~Math.~Phys.\ {\bf 15} (1971) 173
\item{[20]} D.~Buchholz, K.~Fredenhagen: {\it Locality and the structure
of particle states}, Commun.~Math.~Phys.\ {\bf 84} (1982) 1
\item{[21]} S.~Doplicher, J.~E.~Roberts:
 {\it Why there is a field
algebra with a compact gauge group describing the the superselection
structure in particle physics}, Commun.~Math.~Phys.\ {\bf 131} (1990) 51
\item{[22]} J.~Fr{\"o}hlich, G.~Morchi, F.~Strocchi: {\it
Charged sectors and scattering states in quantum electrodynamics},
Ann.~Phys.~(N.Y.) {\bf 119} (1979) 241
\item{[23]} S.~W.~Hawking, G.~F.~R.~Ellis: {\it The Large Scale Structure
of Space-Time}, Cambridge University Press, Cambridge (1973)
\item{[24]} B.~O'Neill: {\it Semi-Riemannian Geometry}, Academic Press,
New York (1983)
\item{[25]} B.~S.~Kay: {\it Linear spin-zero fields in external 
gravitational and scalar fields II}, Commun.~Math.~Phys.\ {\bf 71} (1980) 29
\item{[26]} J.~Dimock: {\it Algebras of local observables on a manifold},
Commun.~Math.~Phys.\ {\bf 77} (1980) 219
\item{[27]} M.~Wollenberg: {\it Scaling limits and type of local
algebras over curved spacetime}. In: W.~B.~Arveson et al.\ (eds.):
{\it Operator Algebras and Topology}, Pitman Research Notes in 
Mathemtics 270, Longman, Harlow (1992)
%%%%%%%%%%%%%%%%%%%%%%%%%%%%%%%%%%%%%%%%%%%%%%%%%%%%%%%%%%%%%%%%%%%%%%%%%%
% In the theory of local observables one starts with a family of von~Neumann
% algebras [1] ${\Cal M}({\Cal O})$ in the vacuum sector.
%
% \proclaim{1.1 Theorem}
% The quotient surface $X$ is  isomorphic to
% a finite union of compact Riemann surfaces
% with a finite number of points removed.
% \endproclaim
%
% \demo{Proof}
% The proof of the theorem.
% \enddemo
% \subheading{ Remarks on the edge of the wedge problem}
% The theory of several complex variables is an important tool
% in quantum field theory.
%
% This calculation implies the following property about real lines:
% \proclaim{1.2 Theorem}
% Let $T^+\subset\Bbb C^2$ ...
% \endproclaim
%%%
% \heading
% 2. Notations and examples
% \endheading
% We start by introducing the main concepts.
% \proclaim{2.1 Definition}
% Let $\Omega$ be as described in the introduction.
% \endproclaim
%
% First we look at the set ${\frak Sem}({\Cal M})$.
% \proclaim{2.2 Theorem}
% The map
% $$ {\frak G}({\Cal M}) \longrightarrow {\frak Sem}({\Cal M})$$
% is a bijection.
% \endproclaim
%%%%%%%%%%%%%
% \heading
% References
% \endheading
%
% \item{[1]} M.~Takesaki: {\it Tomita's Theory of Modular Hilbert Algebras
% and its Applications}, Lecture Notes in Mathematics, Vol.~{\bf 128}
% Springer-Verlag Berlin, Heidelberg, New York (1970).
%
% \item{[2]} E.~Witten: {\it  Topological sigma models},
% Comm.~Math.~Phys., {\bf  118} (1988), 201-210.
%
 \enddocument
 \bye